\documentclass[11pt]{article}
\usepackage{amsfonts}
\usepackage{mathrsfs}
\usepackage{amsmath}
\usepackage{amssymb}
\usepackage[mathscr]{eucal}
\usepackage{graphicx}
\renewcommand{\paragraph}{\roman{paragraph}}

\def \ov{\overline}

\newtheorem{theorem}{\scshape \mdseries  Theorem}[section]
\newtheorem{lemma}[theorem]{\scshape \mdseries  Lemma}
\newtheorem{coro}[theorem]{\scshape \mdseries  Corollary}

\begin{document}

\title{\sf Spectral conditions of complement for some graphical properties }
\author{Guidong Yu\thanks{Email: guidongy@163.com.
Supported by the Natural Science Foundation of Department of
Education of Anhui
Province of China under Grant nos. KJ2015ZD27, KJ2017A362.}, Yi Fang, Yi Xu \\
  {\small  \it School of Mathematics \& Computation Sciences, Anqing Normal University, Anqing 246133,
  China}}

\date{}
\maketitle

\noindent {\bf Abstract:} L.H. Feng at el \cite{feng4} present
sufficient conditions based on spectral radius for a graph with
large minimum degree to be $s$-path-coverable and $s$-Hamiltonian.
Motivated by this study, in this paper, we give the sufficient
conditions for a graph with large minimum degree to be
$s$-connected, $s$-edge-connected, $\beta$-deficient,
$s$-path-coverable, $s$-Hamiltonian and $s$-edge-Hamiltonian in
terms of spectral radius of its complement.

\noindent {\bf Keywords:} Spectral radius; Minimum degree;
Complement; Stability

\noindent {\bf MR Subject Classifications:}  05C50, 05C45, 05C35

\section{Introduction}
Let $G=(V,E)$ be a simple connected graph of order $n$ with vertex
set $V=V(G)=\{v_1,v_2,\ldots,v_n\}$ and edge set $E=E(G)$. The
complement of $G$ is denoted by $\ov{G}$. A regular graph is one
whose vertices all have the same degrees, and semi-regular bipartite
graph is a bipartite graph for which the vertices in the same part
have the same degrees. Let $K_n$, $O_n$ denote the complete graph,
the empty graph on $n$ vertices, respectively. For two disjoint
graphs $G_1$ and $G_2$, the union of $G_1$ and $G_2$, denoted by
$G_1+G_2$, is defined as $V(G_1+G_2)=V(G_1)\cup V(G_2)$ and
$E(G_1+G_2)=E(G_1)\cup(G_2)$; and the join of $G_1$ and $G_2$,
denoted by $G_1 \vee G_2$, is defined as $V(G_1\vee G_2)=V(G_1)\cup
V(G_2)$, and $E(G_1\vee G_2)=E(G_1+G_2)\cup\{xy:x\in V(G_1),y\in
V(G_2)\}$. Denote $K_{n,m}=O_n\vee O_m$, a complete bipartite graph.

The {\it adjacency matrix} of $G$ is defined to be a matrix
$A(G)=[a_{ij}]$ of order $n$, where $a_{ij}=1$ if $v_{i}$ is
adjacent to $v_{j}$, and $a_{ij}=0$ otherwise. The largest
eigenvalue of $A(G)$, denoted by $\mu(G)$, is the spectral radius of
$A(G)$. The degree matrix of $G$ is denoted by
$D(G)=\hbox{diag}\left(d_G(v_1),d_G(v_2), \ldots, d_G(v_n)\right)$,
where $d_G(v)$ or simply $d(v)$ denotes  the degree of a vertex $v$
in the graph $G$, the minimum
   degree of $G$ is denoted by $\delta(G)$. The matrix $L(G)=D(G)-A(G)$ is the {\it
Laplacian matrix} of $G$, and the matrix $Q(G)=D(G)+A(G)$ is the
{\it signless Laplacian matrix} (or {\it $Q$-matrix}) of $G$.

A graph $G$ is said to be {\it $s$-connected} if it has more than
$s$ vertices and remains connected whenever fewer than $s$ vertices
are deleted.

A graph $G$ is {\it $s$-edge-connected} if it has at least two
vertices and remains connected whenever fewer than $s$ edges are
deleted.

The deficiency of a graph $G$, denoted by def($G$), is the number of
vertices unmatched under a maximum matching in $G$. In particular,
$G$ has a 1-factor if and only if def($G$)=0. We say that $G$ is
$\beta$-deficient if def($G$)$\leq\beta$.

A graph $G$ is {\it $s$-path-coverable} if $V(G)$ can be covered by
$s$ or fewer vertex-disjoint paths. In particular, a graph $G$ is
1-path-coverable is the same as $G$ is traceable.

A graph $G$ is {\it $s$-Hamiltonian} if for all $X\subset V(G)$ with
$|X|\leq s$, the subgraph induced by $V(G)\backslash X$ is
Hamiltonian. Thus a graph $G$ is 0-Hamiltonian is the same as $G$ is
Hamiltonian.

A graph $G$ is {\it $s$-edge-Hamiltonian} if any collection of
vertex-disjoint paths with at most $s$ edges altogether belong to a
Hamiltonian cycle in $G$.

The study of the relationship between graph properties and
eigenvalues has attracted much attention. This is largely due to the
following problem of Brualdi and Solheid \cite{Brualdi}: Given a set
$\varrho$ of graphs, find an upper bound for the spectral radii of
the graphs of $\varrho$, and characterize the graphs for which the
maximal spectral radius is attained.  An important paper relating
connectivity and Laplacian eigenvalues of graphs is by Fiedler
\cite{fiedler}, and since then, this relationship has been well
studied. Yu et al \cite{yu3} presented spectral condition for a
graph to be $k$-connected. Cioab・a and Gu \cite{S.M.3} obtained a
sufficient condition for a connected graph to be k-edge-connected in
relation to the second largest eigenvalue. Other related results can
be found in \cite{S.M., gu, liu, s.o}. The study of eigenvalues and
the matching number was initiated by Brouwer and Haemers
\cite{Brouwer}, and then subsequently developed by Cioab・a and many
other researchers \cite{S.M., S.M.1}. Fiedler and Nikiforov
\cite{nikiforov} were the first to establish the spectral sufficient
conditions for graphs to be Hamiltonian or traceable. Zhou
\cite{zhou} presented the signless Laplacian spectral sufficient
conditions for graphs to be Hamiltonian or traceable. Li \cite{li2}
reported several spectral sufficient conditions for Hamiltonian
properties of graphs. Since then, many researchers have studied the
analogous problems under various other conditions; see \cite{Ning0,
Ning1, ning, yu1, feng3, feng4, liu1, lu, nikiforov1, yu,
yu4, yu5, yu2}. Particularly, L.H. Feng at el \cite{feng4} present
sufficient conditions based on spectral radius for a graph with
large minimum degree to be $s$-path-coverable and $s$-Hamiltonian.
Motivated by this study, in this paper, we give the sufficient
conditions for a graph with large minimum degree to be
$s$-connected, $s$-edge-connected, $\beta$-deficient,
$s$-path-coverable, $s$-Hamiltonian and $s$-edge-Hamiltonian in
terms of spectral radius of its complement.

\section{Preliminaries}
For an integer $k\geq0$, the $k$-closure of a graph $G$, denoted by
$C_k(G)$, is the graph obtained from $G$ by successively joining
pairs of nonadjacent vertices whose degree sum is at least $k$ until
no such pair remains,see \cite{bondy}. The $k$-closure of the graph
$G$ is unique, independent of the order in which edges are added.
Note that $d_{C_k(G)}(u)+d_{C_k(G)}(v)\leq k-1$ for any pair of
nonadjacent vertices $u$ and $v$ of $C_k(G)$.

\begin{lemma} \cite{bondy} Let $G$ be a graph of order $n$. Then

(i) $G$ is $s$-connected if and only if $C_{n+s-2}(G)$ is so.

(ii) $G$ is $s$-edge-connected if and only if $C_{n+s-2}(G)$ is so.

(iii) $G$ is $\beta$-deficient if and only if $C_{n-\beta-1}(G)$ is
so.

(iv) $G$ is $s$-path-coverable if and only if $C_{n-s}(G)$ is so.

(v) $G$ is $s$-hamiltonian if and only if $C_{n+s}(G)$ is so.

(vi) $G$ is $s$-edge-hamiltonian if and only if $C_{n+s}(G)$ is so.

\end{lemma}

\begin{lemma} \cite{ning}
Let $G$ be a graph with non-empty edge set. Then $$\mu(G)\geq
min\{\sqrt{d(u)d(v)}:uv\in E(G)\}. $$ Moreover, if $G$ is connected,
then equality holds if and only if $G$ is regular or semi-regular
bipartite.
\end{lemma}

\section{The main results}

Let $EP_n$ be the set of following graphs of order $n$:

$G_1\vee G_2$, where $G_1$ with order $n-k+r$ is a $r$ regular
graph, $G_2$ is a spanning subgraph of $K_{k-r}$.

Let $EC_n$ be the set of following graphs of order $n$:

 $\ov G_1\vee G_2$, where $G_1=(X,Y)$ with order $n-s+1+r$ is a
 semi-regular bipartite graph with $\forall u\in X,~d_{G_{1}}(u)=k-s+2,~\forall v\in Y,~d_{G_{1}}(v)=n-k-1$, $|V(G_2)|=s-1-r\geq 0$.

Let $ES_n$ be the set of following graphs of order $n$:

$\ov G_1\vee G_2$, where $G_1=(X,Y)$ with order $n-s-1+r$ is a
 semi-regular bipartite graph with $\forall u\in X,~d_{G_{1}}(u)=k-s,~\forall v\in Y,~d_{G_{1}}(v)=n-k-1$, $|V(G_2)|=s+1-r\geq 0$.

\begin{theorem}
Let $s\geq1$, $k\geq 1$, $k-s+1\geq0$, $n\geq2k+1$. Let $G$ be a
connected graph of order $n$ and minimum degree $\delta(G)\geq k$.
If
$$\mu(\ov{G})\leq\sqrt{(k-s+2)(n-k-1)},$$ then $G$ is $s$-connected unless
$G\in EP_n$ or $G\in EC_n$.
\end{theorem}

{\bf Proof.} Let $H=C_{n+s-2}(G)$. If $H$ is $s$-connected, then so
is $G$ by Lemma 2.1(i). Now we assume that $H$ is not $s$-connected.
Note that $H$ is $(n+s-2)$-closed, thus every two nonadjacent
vertices $u$, $v$ have degree sum at most $n+s-3$. i.e., for any
edge $uv\in E(\ov{H})$
$$d_{\ov{H}}(u)+d_{\ov{H}}(v)\geq 2(n-1)-(n+s-3)= n-s+1,\eqno(3.1)$$
this implies that
$$d_{\ov{H}}(u)d_{\ov{H}}(v)\geq
d_{\ov{H}}(u)(n-s+1-d_{\ov{H}}(u)).$$

Since $d_H(u)\geq d_G(u)\geq k$ and $d_H(v)\geq d_G(v)\geq k$, we
have $d_{\ov{H}}(u)\leq n-k-1$ and $d_{\ov{H}}(v)\leq n-k-1$.
Combining (3.1), thus $d_{\ov{H}}(u)\geq n-s+1-d_{\ov{H}}(v)=k-s+2$
and $d_{\ov{H}}(v)\geq k-s+2$.

Let $f(x)=x((n-s+1)-x)$ with $k-s+2\leq x\leq n-k-1$, we note $f(x)$
is convex in $x$. We have $f(x)\geq f(k-s+2)$ (or $f(n-k-1)$),
namely $f(x)\geq (k-s+2)(n-k-1)$, the equality holds if and only if
$x=k-s+2$ (or $x=n-k-1$). So $d_{\ov{H}}(u)d_{\ov{H}}(v)\geq
d_{\ov{H}}(u)(n-s+1-d_{\ov{H}}(u))\geq (k-s+2)(n-k-1)$ with equality
if and only if $d_{\ov{H}}(u)=k-s+2$ and $d_{\ov{H}}(v)=n-k-1$. By
Lemma 2.2, Perron-Frobenius theorem, and the assumption,
$$\sqrt{(k-s+2)(n-k-1)}\geq\mu(\ov{G})\geq\mu(\ov{H})\geq\min_{uv\in
E(\ov{H})}\sqrt{d_{\ov{H}}(u)d_{\ov{H}}(v)}\geq\sqrt{(k-s+2)(n-k-1)}.$$
Therefore, $\mu(\ov{G})=\mu(\ov{H})=\sqrt{(k-s+2)(n-k-1)}$, and
$d_{\ov{H}}(u)=k-s+2$, $d_{\ov{H}}(v)=n-k-1$ for any edge $uv\in
E(\ov{H})$. Note that every non-trivial component of $\ov{H}$ has a
vertex of degree at least $n-k-1$ and hence of order at least $n-k$.
This implies that $\ov{H}$ has exactly one nontrivial component $F$
by $n\geq2k+1$. By lemma 2.2, $F$ is either regular or semi-regular
bipartite, where $n-k\leq|V(F)|\leq n$.

Suppose that $F$ is a semi-regular bipartite graph. Then
$n-s+1\leq|V(F)|\leq n$, $F=(X,Y)$ is a semi-regular bipartite graph
where $|X|=n-k-1+m,|Y|=k-s+t+2$ and $0\leq m+t\leq s-1$. For any
vertex $u\in X$, $d_F(u)=k-s+2$, for any vertex $v\in Y$,
$d_F(v)=n-k-1$. We have $\ov H=F+\ov K_{s-1-m-t}$, since
$\mu(\ov{G})=\mu(\ov{H})$ and $\ov{H}$ is a spanning subgraph of
$\ov{G}$, we have $F+\ov K_{s-1-m-t}\subseteq \ov{G}\subseteq F+
K_{s-1-m-t}$, and then $G\in\ov{F}\vee G_1$, where $G_1$ is the
spanning subgraph of $K_{s-1-m-t}$, i.e., $G\in EC_n$; a
contradiction.

Finally we assume $F$ is regular of degree $k-s+2=n-k-1$ when $n=
2k-s+3$. Let $|V(F)|=n-k+r, 0\leq r\leq k$, $\ov{H}=F+\ov{K_{k-r}}$.
Thus $H=\ov F\vee K_{k-r}$, so $\forall u\in\ov F, d_{H}(u)=k,
\forall v\in K_{k-r}, d_{H}(v)=n-1$ in $H$. Since $G\subseteq H,
\delta(G)\geq k$, we get $G\in\ov{F}\vee G_1$, where $G_1$ is the
spanning subgraph of $K_{k-r}$, i.e., $G\in EP_n$.\hfill
$\blacksquare$

Next, we consider the $s$-edge-connected property. Since every
$s$-connected graph is also  $s$-edge-connected. By Theorem 3.1. We
get the following corollary.

\begin{coro}
Let $s\geq1$, $k\geq 1$, $k-s+1\geq0$, $n\geq2k+1$. Let $G$ be a
connected graph of order $n$ and minimum degree $\delta(G)\geq k$.
If
$$\mu(\ov{G})\leq\sqrt{(k-s+2)(n-k-1)},$$ then $G$ is $s$-edge-connected unless
$G\in EP_n$ or $G\in EC_n$.
\end{coro}

\begin{theorem}
Let $k\geq1$, $k\geq2\beta$, $n\geq2k+\beta+2$ with
$n\equiv\beta$(mod2), $0\leq\beta\leq n$. Let $G$ be a connected
graph of order $n$ and minimum degree $\delta(G)\geq k$. If
$$\mu(\ov{G})\leq\sqrt{(\beta+k+1)(n-k-1)},$$ then $G$ is $\beta$-deficient unless
$G\in EP_n$ or $G=K_{k+1}+K_{n-k-1}$ when $\beta=0$.
\end{theorem}

{\bf Proof.} Let $H=C_{n-\beta-1}(G)$. If $H$ is $\beta$-deficient,
then so is $G$ by Lemma 2.1(iii). Now we assume that $H$ is not
$\beta$-deficient. Note that $H$ is $(n-\beta-1)$-closed, thus every
two nonadjacent vertices $u$, $v$ have degree sum at most
$n-\beta-2$, i.e.,
$$d_{\ov{H}}(u)+d_{\ov{H}}(v)\geq 2(n-1)-(n-\beta-2)= n+\beta,\eqno(3.2)$$
for any edge $uv\in E(\ov{H})$.

Since $d_H(u)\geq d_G(u)\geq k$ and $d_H(v)\geq d_G(v)\geq k$, we
have $d_{\ov{H}}(u)\leq n-k-1$ and $d_{\ov{H}}(v)\leq n-k-1$.
Combining (3.2), thus $d_{\ov{H}}(u)\geq
n+\beta-d_{\ov{H}}(v)=\beta+k+1$ and $d_{\ov{H}}(v)\geq \beta+k+1$.
This implies that
$$d_{\ov{H}}(u)d_{\ov{H}}(v)\geq
d_{\ov{H}}(u)(n+\beta-d_{\ov{H}}(u)).$$ Let $f(x)=x((n+\beta)-x)$
with $\beta+k+1\leq x\leq n-k-1$, we note $f(x)$ is convex in $x$.
We have $f(x)\geq f(\beta+k+1)$ (or $f(n-k-1)$), namely $f(x)\geq
(\beta+k+1)(n-k-1)$, the equality holds if and only if $x=\beta+k+1$
(or $x=n-k-1$). So $d_{\ov{H}}(u)d_{\ov{H}}(v)\geq
d_{\ov{H}}(u)(n+\beta-d_{\ov{H}}(u))\geq (\beta+k+1)(n-k-1)$ with
equality if and only if $d_{\ov{H}}(u)=\beta+k+1$ and
$d_{\ov{H}}(v)=n-k-1$. By Lemma 2.2, Perron-Frobenius theorem, and
the assumption,
$$\sqrt{(\beta+k+1)(n-k-1)}\geq\mu(\ov{G})\geq\mu(\ov{H})\geq\min_{uv\in
E(\ov{H})}\sqrt{d_{\ov{H}}(u)d_{\ov{H}}(v)}\geq\sqrt{(\beta+k+1)(n-k-1)}.$$
Therefore, $\mu(\ov{G})=\mu(\ov{H})=\sqrt{(\beta+k+1)(n-k-1)}$, and
$d_{\ov{H}}(u)=\beta+k+1$, $d_{\ov{H}}(v)=n-k-1$  for any edge
$uv\in E(\ov{H})$. Note that every non-trivial component of $\ov{H}$
has a vertex of degree at least $n-k-1$ and hence of order at least
$n-k$, this implies that $\ov{H}$ has exactly one nontrivial
component $F$ by $n\geq2k+\beta+2$. By lemma 2.2, $F$ is either
regular or semi-regular bipartite, where $n-k\leq|V(F)|\leq n$.

Suppose that $F$ is a semi-regular bipartite graph. Then
$n+\beta\leq|V(F)|\leq n$, so $\beta=0$ and
$F=\ov{H}=K_{k+1,n-k-1}$. Since $\mu(\ov{G})=\mu(\ov{H})$ and
$\ov{H}$ is a connected spanning subgraph of $\ov{G}$, we have
$\ov{G}=\ov{H}$ by Perron-Frobenius theorem. Thus
$\ov{G}=K_{k+1,n-k-1}$, this implies that $G=K_{k+1}+K_{n-k-1}$; a
contradiction.

Finally we assume $F$ is regular of degree $\beta+k+1=n-k-1$ when
$n= 2k+\beta+2$, so $n-k\leq|V(F)|\leq n$. Let $|V(F)|=n-k+r, 0\leq
r\leq k$, $\ov{H}=F+\ov {K_{k-r}}$. Thus $H=\ov F\vee K_{k-r}$, so
$\forall u\in\ov F, d_{H}(u)=k, \forall v\in K_{k-r}, d_{H}(v)=n-1$
in $H$. Since $G\subseteq H, \delta(G)\geq k$, we get
$G\in\ov{F}\vee G_1$, where $G_1$ is the spanning subgraph of
$K_{k-r}$, $G\in EP_n$; a contradiction. \hfill $\blacksquare$

\begin{theorem}
Let $s\geq1$, $k\geq1$, $n\geq2k+s+1$. Let $G$ be a connected graph
of order $n$ and minimum degree $\delta(G)\geq k$. If
$$\mu(\ov{G})\leq\sqrt{(k+s)(n-k-1)},$$ then $G$ is $s$-path-coverable
unless $G\in EP_n$ or $G=K_{k+1}+K_{n-k-1}$ when $s=1$.
\end{theorem}

 {\bf Proof.} Let $H=C_{n-s}(G)$. If $H$ is $s$-path-coverable, then so is
$G$ by Lemma 2.1(iv). Now we assume that $H$ is not
$s$-path-coverable. Note that $H$ is $(n-s)$-closed, thus every two
nonadjacent vertices $u$, $v$ in $H$ have degree sum at most
$n-s-1$, i.e., for any edge $uv\in E(\ov{H})$
$$d_{\ov{H}}(u)+d_{\ov{H}}(v)\geq 2(n-1)-(n-s-1)= n+s-1,\eqno(3.3)$$
this implies that
$$d_{\ov{H}}(u)d_{\ov{H}}(v)\geq
d_{\ov{H}}(u)(n+s-1-d_{\ov{H}}(u)).$$

Since $d_{H}(u)\geq d_{G}(u)\geq k$ and $d_{H}(v)\geq d_{G}(v)\geq
k$, we have $d_{\ov{H}}(u)\leq n-k-1$ and $d_{\ov{H}}(v)\leq n-k-1$.
Combining (3.3), thus $d_{\ov{H}}(u)\geq n+s-1-d_{\ov{H}}(u)=k+s$
and $d_{\ov{H}}(v)\geq k+s$.

Let $f(x)=x((n+s-1)-x)$, $k+s\leq x\leq n-k-1$. We note $f(x)$ is
convex in $x$. Then we have $f(x)\geq f(k+s)$ (or $f(n-k-1)$),
namely $f(x)\geq (k+s)(n-k-1)$, the equality holds if and only if
$x=k+s$ (or $x=n-k-1$). So $d_{\ov{H}}(u)d_{\ov{H}}(v)\geq
d_{\ov{H}}(u)(n+s-1-d_{\ov{H}}(u))\geq (k+s)(n-k-1)$ with equality
if and only if $d_{\ov{H}}(u)=k+s$ and $d_{\ov{H}}(v)=n-k-1$. By
Lemma 2.2, Perron-Frobenius theorem, and the assumption,
$$\sqrt{(k+s)(n-k-1)}\geq\mu(\ov{G})\geq\mu(\ov{H})\geq\min_{uv\in
E(\ov{H})}\sqrt{d_{\ov{H}}(u)d_{\ov{H}}(v)}\geq\sqrt{(k+s)(n-k-1)}.$$
Therefore, $\mu(\ov{G})=\mu(\ov{H})=\sqrt{(k+s)(n-k-1)}$, and
$d_{\ov{H}}(u)=k+s$, $d_{\ov{H}}(v)=n-k-1$ for any edge $uv\in
E(\ov{H})$. Note that every non-trivial component of $\ov{H}$ has a
vertex of degree at least $n-k-1$ and hence of order at least $n-k$,
this implies that $\ov{H}$ has exactly one nontrivial connected
component $F$ by $n\geq 2k+s+1$. By lemma 2.2, $F$ is either regular
or semi-regular bipartite, where $n-k\leq|V(F)|\leq n$.

Suppose that $F$ is a semi-regular bipartite graph. Then
$n+s-1\leq|V(F)|\leq n$, and then $s=1$, $F=\ov{H}$. Since
$\mu(\ov{G})=\mu(\ov{H})$ and $\ov{H}$ is a connected spanning
subgraph of $\ov{G}$, we have $\ov{G}=\ov{H}$ by Perron-Frobenius
theorem. Thus $\ov{G}=K_{k+1,n-k-1}$, this implies that
$G=K_{k+1}+K_{n-k-1}$; a contradiction.

Next suppose $F$ is regular of degree $k+s=n-k-1$, then $n=2k+s+1$.
Let $|V(F)|=n-k+r,~0\leq r\leq k$, $\ov{H}=F+\ov{K_{k-r}}$. Thus
$H=\ov F\vee K_{k-r}$. Because  $\forall u\in\ov F, d_{H}(u)=k,
\forall v\in K_{k-r}, d_{H}(v)=n-1$ in $H$, $G\subseteq H$ and
$\delta(G)\geq k$, we get $G=\ov{F}\vee G_2$, where $G_2$ is the
spanning subgraph of $K_{k-r}$, and then $G\in EP_n$; a
contradiction.\hfill $\blacksquare$

\begin{theorem}
Let $s\geq0$, $k\geq s+1$, $n\geq2k+1$. Let $G$ be a connected graph
of order $n$ and minimum degree $\delta(G)\geq k$. If
$$\mu(\ov{G})\leq\sqrt{(k-s)(n-k-1)}.$$ Then $G$ is $s$-hamiltonian unless
$G\in EP_n$ or $G\in ES_n$.
\end{theorem}

{\bf Proof.} Let $H=C_{n+s}(G)$. If $H$ is $s$-hamiltonian, then so
is $G$ by Lemma 2.1(v). Now we assume that $H$ is not
$s$-hamiltonian. Note that $H$ is $(n+s)$-closed, thus every two
nonadjacent vertices $u$, $v$ have degree sum at most $n+s-1$, i.e.,
$$d_{\ov{H}}(u)+d_{\ov{H}}(v)\geq 2(n-1)-(n+s-1)= n-s-1,\eqno(3.4)$$
for any edge $uv\in E(\ov{H})$.  This implies that
$$d_{\ov{H}}(u)d_{\ov{H}}(v)\geq
d_{\ov{H}}(u)(n-s-1-d_{\ov{H}}(u)).$$

Since $d_H(u)\geq d_G(u)\geq k$ and $d_H(v)\geq d_G(v)\geq k$, we
have $d_{\ov{H}}(u)\leq n-k-1$ and $d_{\ov{H}}(v)\leq n-k-1$.
Combining (3.4), thus $d_{\ov{H}}(u)\geq n-s-1-d_{\ov{H}}(v)=k-s$
and $d_{\ov{H}}(v)\geq k-s$. Let $f(x)=x((n-s-1)-x)$ with $k-s\leq
x\leq n-k-1$, we note $f(x)$ is convex in $x$. We have $f(x)\geq
f(k-s)$ (or $f(n-k-1)$), namely $f(x)\geq (k-s)(n-k-1)$, the
equality holds if and only if $x=k-s$ (or $x=n-k-1$). So
$d_{\ov{H}}(u)d_{\ov{H}}(v)\geq
d_{\ov{H}}(u)(n-s-1-d_{\ov{H}}(u))\geq (k-s)(n-k-1)$ with equality
if and only if $d_{\ov{H}}(u)=k-s$ and $d_{\ov{H}}(v)=n-k-1$. By
Lemma 2.2, Perron-Frobenius theorem, and the assumption,
$$\sqrt{(k-s)(n-k-1)}\geq\mu(\ov{G})\geq\mu(\ov{H})\geq\min_{uv\in
E(\ov{H})}\sqrt{d_{\ov{H}}(u)d_{\ov{H}}(v)}\geq\sqrt{(k-s)(n-k-1)}.$$
Therefore, $\mu(\ov{G})=\mu(\ov{H})=\sqrt{(k-s)(n-k-1)}$, and
$d_{\ov{H}}(u)=k-s$, $d_{\ov{H}}(v)=n-k-1$ for any edge $uv\in
E(\ov{H})$. Note that every non-trivial component of $\ov{H}$ has a
vertex of degree at least $n-k-1$ and hence of order at least $n-k$,
This implies that $\ov{H}$ has exactly one nontrivial component $F$
by $n\geq2k+1$. By lemma 2.2,  $F$ is either regular or semi-regular
bipartite, where $n-k\leq|V(F)|\leq n$.

Suppose that $F$ is a semi-regular bipartite graph. Then
$n-s-1\leq|V(F)|\leq n$, $F=(X,Y)$ is a semi-regular bipartite graph
where $|X|=n-k-1+m,|Y|=k-s+t$ and $0\leq m+t\leq s+1$. For any
vertex $u\in X$, $d_F(u)=k-s$, for any vertex $v\in Y$,
$d_F(v)=n-k-1$. We have $\ov H=F+\ov K_{s+1-m-t}$, since
$\mu(\ov{G})=\mu(\ov{H})$ and $\ov{H}$ is a spanning subgraph of
$\ov{G}$, we have $F+\ov K_{s+1-m-t}\subseteq \ov{G}\subseteq F+
K_{s+1-m-t}$, and then $G\in\ov{F}\vee G_1$, where $G_1$ is the
spanning subgraph of $K_{s+1-m-t}$, i.e., $G\in ES_n$; a
contradiction.

Finally we assume $F$ is regular of degree $k-s=n-k-1$ when $n=
2k-s+1$.  Let $|V(F)|=n-k+r, 0\leq r\leq k$,
$\ov{H}=F+\ov{K_{k-r}}$. Thus $H=\ov F\vee K_{k-r}$, so $\forall
u\in\ov F, d_{H}(u)=k, \forall v\in K_{k-r}, d_{H}(v)=n-1$ in $H$.
Since $G\subseteq H, \delta(G)\geq k$, we get $G\in\ov{F}\vee G_1$,
where $G_1$ is the spanning subgraph of $K_{k-r}$, i.e., $G\in
EP_n$.\hfill $\blacksquare$

Next, we consider the $s$-edge-hamiltonian property. Since every
$s$-hamiltonian graph is also  $s$-edge-hamiltonian. By Theorem 3.5,
we get the following corollary.

\begin{coro}
Let $s\geq0$, $k\geq s+1$, $n\geq2k+1$. Let $G$ be a connected graph
of order $n$ and minimum degree $\delta(G)\geq k$. If
$$\mu(\ov{G})\leq\sqrt{(k-s)(n-k-1)},$$ then $G$ is $s$-edge-hamiltonian unless
$G\in EP_n$ or $G\in ES_n$.
\end{coro}


\vspace{3mm}
 {\small \begin{thebibliography}{90}

 \bibitem{bondy} J.A. Bondy, V. Chvatal, A method in graph theory,
{\it Discrete Math.}, 15(2) (1976), 111-135.

\bibitem{Brualdi}  R.A. Brualdi, E.S. Solheid, On the spectral radius of complementary acyclic matrices of zeros and
ones, {\it SIAM J. Algebr. Discrete Methods}, 7 (1986), 265-272.

\bibitem{Brouwer} A. Brouwer, W. Haemers, Eigenvalues and perfect matchings,{\it Linear Algebra
Appl.}, 395 (2005), 155-162.

\bibitem{zhang} A. Berman, X.D. Zhang, On the spectral radius of graphs with cut vertices. {\it J. Combin.
Theory Ser. B}, 83 (2001), 233-240.

\bibitem{S.M.} S.M. Cioab・a, Eigenvalues and edge-connectivity of regular graphs, {\it Linear Algebra Appl}.
432 (2010), 458-470.

\bibitem{S.M.1} S.M. Cioab・a, D. Gregory, Large matchings from eigenvalues, {\it Linear Algebra Appl}. 422
(2007), 308-317.

\bibitem{S.M.3} S.M. Cioab・a, X.F. Gu, Connectivity, toughness, spanning trees of bounded degree, and
the spectrum of regular graphs, {\it Czech. Math. Journal}, 2016, 66(3), 913-924.

\bibitem{ning} B.-L. Li, B. Ning, Spectral analogues of
Erd\"{o}s' and Moon-Moser's theorems on Hamilton cycles, {\it Linear
and Multilinear Algebra}, 64 (2016), 2252-2269.

\bibitem{liu} H. Liu, M. Lu, F. Tian, Edge-connectivity and (signless) Laplacian eigenvalue of graphs,
{\it Linear Algebra Appl.}, 439 (2013), 3777-3784.

\bibitem{lu} M. Lu, H. Liu, F. Tian, Spectral radius and Hamiltonian graphs, {\it Linear Algebra
Appl.}, 437 (2012), 1670-1674.

\bibitem{liu1}  R. Liu, W.C. Shiu, J. Xue, Sufficient spectral conditions on Hamiltonian and traceable
graphs, {\it Linear Algebra Appl.}, 467 (2015), 254-266.

\bibitem{li2} R. Li, Eigenvalues, Laplacian eigenvalues and some Hamiltonian properties of graphs,
{\it Util. Math.}, 88 (2012), 247-257.

\bibitem{Ning0} B. Ning, J. Ge, Spectral radius and Hamiltonian properties of
graphs, {\it Linear and Multilinear Algebra}, 63(8)(2015),
1520-1530.

\bibitem{Ning1} B. Ning, B. Li, Spectral radius and traceability of connected claw-free graphs, {\it Mathematics} 2(2014)(2015), 1-8.

\bibitem{nikiforov1} V. Nikiforov, Spectral radius and Hamiltonicity
of graphs with large minimum degree, {\it Czechoslovak Mathematical
Journal}, 66(3)(2016), 925-940.

\bibitem{feng3} L.H. Feng, P.L. Zhang, H. Liu, W.J. Liu, M.M. Liu, Y.Q. Hu, Spectral conditions for some graphical properties,
{\it Linear Algebra Appl.}, 524 (2017), 182-198.

\bibitem{feng4} L.H. Feng, W.J. Liu, M.M. Liu, P.L. Zhang, Spectral conditions for graphs to be $k$-Hamiltonian or $k$-path-coverable,
ResearchGate.

\bibitem{fiedler} M. Fiedler, Algebraic connectivity of graphs, {\it Czechoslovak Math. J.}, 23 (1973), 298-305.

\bibitem{nikiforov} M. Fiedler, V. Nikiforov, Spectral radius and Hamiltonicity of graphs, {\it Linear Algebra
Appl.}, 432 (2010), 2170-2173.

\bibitem{yu1} Y.-Z. Fan, G.-D. Yu, Spectral condition for a graph to be Hamiltonian with respect to normalized Laplacian,
{\it Mathematics}, (2012).

\bibitem{gu} X. Gu, H.-J. Lai, P. Li, S. Yao, Edge-disjoint spanning trees, edge connectivity and
eigenvalues in graphs, J. Graph Theory 81 (2016), 16-29.

\bibitem{s.o} O. Suil, S.M. Cioab\v{a}, Edge-connectivity, eigenvalues, and matchings in regular graphs,
{\it SIAM J. Discrete Math.}, 24 (2010), 1470-1481.

\bibitem{yu}G.-D. Yu, Y.-Z. Fan, Spectral conditions for a graph to be
Hamilton-connected, {\it Applied Mechanics and Materials}, 336-338
(2013), 2329-2334.

\bibitem{yu4}G.-D. Yu, M.-L. Ye, G.-X. Cai, J.-D. Cao, Signless Laplacian Spectral
Conditions for Hamiltonicity of Graphs, {\it Journal of Applied
Mathematics}, vol. 2014, Article ID 282053, 6 pages, 2014.

\bibitem{yu5} G.-D. Yu, G.-X. Cai, M.-L. Ye, J.-D. Cao, Energy conditions for Hamiltonicity of graphs, {\it Discrete Dynamics in Nature and
Society}, 2014, Article ID 305164,6 pages.

\bibitem{yu2}G.-D. Yu, Spectral Radius and Hamiltonicity of a Graph, {\it
Mathematic Applicata}, 27(3)(2014), 588-595.

\bibitem{yu3}G.-D. Yu, R. Li, B.-H. Xing, Spectral
Invariants and Some Stable Properties of a Graph, {\it Ars
Combinatoria}, 121(2015), 33-46.

\bibitem{zhou} B. Zhou, Signless Laplacian spectral radius and Hamiltonicity, {\it Linear Algebra Appl.},
432(2010), 566-570.
\end{thebibliography}}
\end{document}